\documentclass[12pt,twoside,a4paper]{amsart}
\usepackage[latin1]{inputenc}
\usepackage[TS1,T1]{fontenc} 
\usepackage{ae}
\usepackage{aeguill}
\usepackage[francais]{babel}
\usepackage{amsfonts,amssymb,latexsym,xspace}  
\usepackage{amsmath,enumerate}  
\usepackage{amsthm}
\usepackage{amscd}
\usepackage{graphicx}
\usepackage{amsfonts,latexsym,epsf}
\numberwithin{equation}{section}
\usepackage{pstricks}
\usepackage{pst-plot}
\usepackage{pst-node}
\def\cal {\mathcal}

\newcommand{\ds}{\displaystyle}
\newcommand{\ntoinf}{n\longrightarrow\infty}
\newcommand{\ItoI}{\mathcal{I}\longrightarrow\mathcal{I}}

\newcommand{\oen}{\overline{E}_{N}}

\newcommand{\wpsi}{\widehat{\psi}}

\newcommand{\ie}{\mathcal{I}}
\newcommand{\W}{\mathcal{W}}
\newcommand{\dji}{\mathcal{G}}
\newcommand{\key}{\mathcal{K}}
\newcommand{\ow}{\mathcal{O}}
\newcommand{\h}{\mathcal{H}}
\newcommand{\Z}{\mathbb{Z}}
\newcommand{\N}{\mathbb{N}}
\newcommand{\C}{\mathbb{C}}
\newcommand{\R}{\mathbb{R}}
\newcommand{\qu}{\mathcal{Q}}

\def \marginpar#1{}
\newtheorem*{abs}{Abstract}
\newtheorem*{rsm}{R{\'e}sum{\'e}}
\newtheorem{prop}{Proposition}

\newtheorem{lem}{Lemme}
\newtheorem{thm}{Th{\'e}or{\`e}me}
\newtheorem{defn}{\bf D{\'e}finition}
\newenvironment{rem}{\noindent {\bf Remarque :} }{}
\def\pn {\par \noindent}
\def\md {\par \medskip}

\title[Th{\'e}or{\`e}me de la limite locale]{\bf Un Th{\'e}or{\`e}me de la limite locale
  pour des Algorithmes Euclidiens}
\author{ A{\"\i}cha {\sc Hachemi}} 

\address{
 A{\"\i}cha {\sc Hachemi} : Laboratoire de G{\'e}om{\'e}trie et Dynamique, 
Institut de Math{\'e}\-ma\-ti\-ques de Jussieu, F-75251 Paris,  France.}

\email{ hachemi@math.jussieu.fr  }
\date{Juin 2004}

\begin{document}
\maketitle

\begin{rsm}
 Nous montrons un th{\'e}or{\`e}me de la limite locale  
pour les
Algorithmes Euclidiens : $standard,$\;$\mbox{centr{\'e}}$ et $impair,$\; avec
une fonction co{\^u}t quelconque {\`a} croissance mod{\'e}r{\'e}e. 
\end{rsm}
\begin{abs}
We prove a local limit theorem 
for the  Euclidean algorithms $standard,$\;$\mbox{centred}$ and
$odd$, with
any cost function of moderate growth.
\end{abs}
\section{Introduction}

\noindent
Il est devenu classique d'associer des syst{\`e}mes dynamiques de l'intervalle
{\`a} des  algorithmes arithm{\'e}tiques (voir B.Vall{\'e}e \cite{Vallee2000}). Dans un article
r{\'e}cent \cite{BVVB2003},
 V.Baladi et B.Vall{\'e}e ont {\'e}tudi{\'e} certaints de ces algorithmes {\`a}
 savoir les algorithmes {\it``standard", ``centr{\'e}",}  et
 {\it ``impair",} associ{\'e}s {\`a} des fonctions co{\^u}t {\`a} croissance mod{\'e}r{\'e}e. 
Gr{\^a}ce {\`a} une
 {\'e}tude du comportement spectral d'un op{\'e}rateur de transfert et {\`a} la
 formule de Perron, elles ont
 obtenu un th{\'e}or{\`e}me de la limite centrale (th{\'e}or{\`e}me 2, ci-dessous) 
avec vitesse de convergence 
optimale, et un th{\'e}or{\`e}me de la limite locale pour des
 fonctions co{\^u}t $``\mbox{\it r{\'e}seau}"$ avec la m{\^e}me vitesse de
 convergence. Nous montrons le deuxi{\`e}me th{\'e}or{\`e}me sans vitesse pour des fonctions
 co{\^u}t quelconques.

\noindent
Les trois algorithmes cit{\'e}s ci-dessus sont d{\'e}finis par des divisions
euclidiennes. 
Soient $u,
v$ deux entiers tels que $v \geq u \geq 1$. La division classique (qui
correspond {\`a} l'algorithme Euclidien standard $\mathcal{G}$), $v=mu+r$
produit un entier $m\geq 1$ et un reste entier  $r$ tel que $0 \leq r <
u$. La division centr{\'e}e (l'algorithme centr{\'e} $\mathcal{K}$) exige que
$v \geq 2u$ et prend la forme $v =mu+s$, avec $s \in [-u/2, +u/2[$. En
posant $s = \varepsilon r$, avec $\varepsilon = \pm 1$ (et $\varepsilon = +1$,
si $s=0$ ), on obtient un entier reste $r$ tel que $0 \leq r \leq u/2$
et un entier $m \geq
2$. La division impaire (l'algorithme $\mathcal{O}$) produit un
quotient impair : $v = mu + s $ avec $m$ impair et $s$ un entier $s
\in [-u, +u[$. En posant $s= \varepsilon r$ avec $\varepsilon = \pm 1$
(et $\varepsilon = +1$, si $s=0$) on obtient un  reste  entier $r$ tel que
$0 \leq r \leq u$ et un entier $m \geq 1$.

\noindent 
Dans les trois cas, les divisions sont d{\'e}finies par des paires $q = (m,
\varepsilon)$, appel{\'e}es $``digits"$.

\noindent Tout couple d'entiers $(u, v)$ engendre une suite de   
``transformations fractionnelles lin{\'e}aires" (TFLs) $h$ dans un
 ensemble $\h$ ($\h$ d{\'e}pend de chaque algorithme), qui transforme le
 quotient $r/u$ en une fonction de
$u/v$. On a $h_{[m, \varepsilon]}(x) = 1/(m + \varepsilon x)$. La TFL
 qui appara{\^\i}t dans la derni{\`e}re {\'e}tape appartient {\`a} $\mathcal{F}\subset\h$
 (l'algorithme s'arr{\^e}te lorsque $r=0$).

\pn 
Ainsi, chaque algorithme appliqu{\'e} {\`a} un rationnel $u/v$ donne une
 fraction continue 
\begin {equation*}
\frac{u}{v}=  \frac{1}{\ds{  m_1+
  \frac{\varepsilon_1}{\ds{  m_2+
  \frac{\varepsilon_2}{~\ds  \ddots +\frac  {\varepsilon_{P-1}} 
{ m_P }~}        }}}} \, ,
\end{equation*} 

\pn de longueur $P = P(u,v)$, et qui d{\'e}compose $u/v$ en 
$$ (u/v) = h_{1} \circ  h_{2} \circ ...\circ h_{P}(0)= h(0)$$
o{\`u} les $h_{i}\in\h, 1\leq i\leq P-1$, et $h_{p}\in\mathcal{F}$.

\pn On s'int{\'e}resse aux diff{\'e}rents co{\^u}ts associ{\'e}s {\`a} l'ex{\'e}cution d'un algorithme. 
Le co{\^u}t le plus basique est le nombre d'{\'e}tapes $P$. En g{\'e}n{\'e}ral, {\'e}tant
donn{\'e}e une fonction co{\^u}t $c$ {\`a} valeurs non-n{\'e}gatives d{\'e}finie sur $\mathcal{H}:=
\mbox{l'ensemble des TFLs}$ associ{\'e}es {\`a} l'algorithme, on consid{\`e}re un co{\^u}t total additif de
la forme
$$C(u, v) : = \sum_{i=1}^{P(u,v)} c(h_i).$$
On associe {\`a} chaque algorithme un syst{\`e}me dynamique de l'intervalle
 $T : \mathcal{I} \longrightarrow \mathcal{I}$ (avec $\mathcal{I}= ]0,
      1[$ pour $\mathcal{G}$,  $\mathcal{I} =]0, 1/2]$ pour
      $\mathcal{K}$ et $\mathcal{I} =[0, 1]$ pour
      $\mathcal{O}$).  $T$ est l'extension {\`a}
      $\mathcal{I}$ de l'application d{\'e}finie sur les rationnels en associant $r/u$ {\`a}  $u/v$.
On obtient
$$T(x) := \mid \frac{1}{x} - A(\frac{1}{x})\mid , \;\;\;\; x \neq
0,\;\;\;T(0)=0$$ 
o{\`u} 
$$\left\{ \begin{array}{l}
A(y):= \mbox{la partie enti{\`e}re de}\; y\; ;\quad \mbox{pour}\; \mathcal{G}\\
A(y):= \mbox{l'entier le plus proche de}\; y\;;\quad \mbox{pour}\; \mathcal{K}\\
A(y):= \mbox{l'entier impair le plus proche de}\; y\;;\quad
\mbox{pour}\; 
\mathcal{O}\\
\end{array} \right.$$  
et l'ensemble $\mathcal{H} = \{ h_{[q]} \}$ est l'ensemble des
branches inverses de $T$. L'ensemble des branches inverses de l'it{\'e}r{\'e}
$T^n$ est $\h^n$, ses {\'e}l{\'e}ments sont de la forme $ h_{[q_1]}\circ
h_{[q_2]}\circ\cdots\circ h_{[q_n]}$ o{\`u} $n$ est appel{\'e} {\it ``profondeur"}
de la branche.
Les syst{\`e}mes $T$ associ{\'e}s aux trois algorithmes $\mathcal{G,
  K,O}$, appartiennent {\`a} la ``bonne classe" des ``applications compl{\`e}tes
par morceaux", que l'on d{\'e}finit comme suit :
\begin{defn}
$\left[\mbox{Application de l'intervalle compl{\`e}te par
    morceaux}\right]$.

\pn Une application $T:\ItoI$ est compl{\`e}te par morceaux s'il existe un
ensemble $\qu$ (fini ou d{\'e}nombrable) 
 et une partition d'ouverts $\{\ie _{q}\}_{q\in\qu}$ (mod un ensemble
d{\'e}nombrable) de l'intervalle $\ie$ tels que la restriction de $T$ {\`a}
$\ie_{q}$ admette une extension bijective de classe $\mathcal{C}^{2}$ de
la cl{\^o}ture de ${\cal I}_{q}$ dans ${\cal I}$.
\end{defn}

\begin{defn}$\left[\mbox{Bonne classe}\right]$.
 Une application compl{\`e}te par morceaux
  appartient {\`a} la bonne classe si:

 $(i)$ $T$ est uniform{\'e}ment dilatante par morceaux, c.{\`a}.d, il
existe $C > 0$ et $\widehat{\rho} < 1$ tels que $\mid h^{'}(x)\mid\leq
C\;\widehat{\rho}^n$\;pour tout $h\in\h^n$, toute profondeur
n et tout $x\in\ie$. Le nombre $\rho$ d{\'e}fini par :
\begin{equation*} 
\rho := \limsup_{\ntoinf}\; (\max\{|h^{'}(x)|;\; h\in\h^{n},
x\in\ie\})^{1/n}
\end{equation*}
est appel{\'e} $\mbox{le taux de contraction}$.

 $(ii)$ Il existe $\widehat {K}>0,$ appel{\'e}e constante 
de distorsion, telle que toute branche
inverse $h$ de $T$ v{\'e}rifie: $$|h^{''}(x)|\leq\widehat {K}|h^{'}(x)|\;\;\;
\mbox{pour tout}\; x\in\ie.$$ 
$(iii)$ Il existe $\sigma_{0} <1$ tel que $\displaystyle \sum_{h\in\h} \sup
|h^{'}|^{\sigma} <\infty$ pour tout r{\'e}el $\sigma >\sigma_{0}$.
 
 $(i\upsilon)$ L'application $T$ n'est pas 
conjugu{\'e}e {\`a} une application affine par morceaux.
\end{defn}

\begin{rem}
On v{\'e}rifie que pour nos algorithmes $\sigma_{0}=1/2$ (voir \cite{Broise}).
\end{rem}

\medskip
\noindent Si $\ie$ est muni d'une probabilit{\'e} (initiale) de densit{\'e} $f_0$ par
rapport {\`a} la mesure de Lebesgue, alors $T$ agit sur $\ie$  par
$(f_{1}\;dx) = (T_{\star}(f_{0}\;dx))$. L'op{\'e}rateur ${\bf H}$ tel que
$f_{1} = {\bf H}[f_{0}]$ est appel{\'e} le transformateur de densit{\'e},
ou l'op{\'e}rateur de Perron-Frobenius. Un changement de variable donne:
$${\bf H}[f](x) := \ds \sum_{h\in\h}|h^{'}(x)|\; f\circ h(x),\;\;\;\;{\bf
  H}^{n}[f](x) :=  \ds \sum_{h\in\h^{n}}|h^{'}(x)|\; f\circ h(x)$$

\medskip
\noindent {\bf Condition $\mathcal{CM}$}[Croissance Mod{\'e}r{\'e}e]. Soit $\h$
l'ensemble des branches inverses d'une application de la bonne
classe. Un co{\^u}t $c : \h\longrightarrow\R^{+}$ est {\`a} croissance
mod{\'e}r{\'e}e si:
\[\ds \sum_{h\in\h}\exp[w c(h)]\;.\;|h^{'}(x)|^s \]
converge lorsque $(\Re s, \Re w)\in\Sigma_{0} \times
W_{0}$\;\;avec\;\;$\Sigma_{0}=]\widehat\sigma_{0}, \infty]$,
\;pour\;$\sigma_{0}\leq\widehat\sigma_{0}<1$, et $W_{0}= ]-\infty,
\nu_{0}]$\;pour\;$\nu_{0} > 0$.

\noindent En  posant $\h^{\star} := \cup_{k\geq 1}\h^{k}$ on peut prolonger le
co{\^u}t en un co{\^u}t total, qui sera aussi not{\'e} $c$, d{\'e}fini sur
$\h^{\star}$ par 
$$c(h_{1}\circ h_{2}\circ\dots\circ h_{k}) := \ds \sum_{i=1}^{k}
c(h_{i}).$$
On peut alors d{\'e}finir une version perturb{\'e}e et pond{\'e}r{\'e}e de l'op{\'e}rateur
de transfert d{\'e}pendant de deux param{\`e}tres complexes $s$\;et\;$w$,

\[{\bf H}_{s,w}[f](x) := \ds \sum_{h\in\h} \exp[w
c(h)]\;.\;|h^{'}(x)|^{s}\;.\; 
f\circ h(x).\]

\pn Par cons{\'e}quent, (en utilisant la propri{\'e}t{\'e} d'additivit{\'e} du
co{\^u}t total $C$)
\[
{\bf H}^{n}_{s,w}[f](x) := \ds \sum_{h\in\h^{n}} 
\exp[w c(h)]\;.\;|h^{'}(x)|^{s}\;.\; f\circ h(x).
\]

\md
\pn Dans la proposition suivante on cite quelques propri{\'e}t{\'e}s de
 l'op{\'e}rateur
 ${\bf H}_{s,w}$.
\begin{prop} \cite{BVVB2003}. En posant $s=\sigma +
  it$, $w= i\tau$, et $R(s, w)$ le rayon spectral de  ${\bf
    H}_{s, w}$ et $R_{e}(s,w)$ son rayon spectral essentiel, et ${\bf
    H}_{\sigma,0}:= {\bf H}_\sigma$, on a :

\md
$(1)$ Si $\sigma\in \Sigma_{0}$, alors ${\bf
  H}_{s,i\tau}$ est born{\'e} sur $C^{1}(\mathcal{I})$, et il d{\'e}pend
analytiquement de $(s,i\tau),\;R(s,i\tau)\leq R(\sigma)$ et
$R_{e}(s,i\tau)\leq \widehat\rho R_{e}(\sigma)$ (avec
$\rho<\widehat\rho<1$). De plus l'op{\'e}rateur ${\bf H}_{\sigma}$
admet une unique valeur propre $\lambda(\sigma)$ r{\'e}elle, positive, simple et
de module maximal, associ{\'e}e {\`a} une fonction propre positive.

$(2)$ {\rm  [Trou Spectral.] }  Pour  $\sigma \in 
\Sigma_0$, il existe un trou spectral, c.{\`a}.d, 
le {\it rayon spectral sous-dominant } 
$r_{\sigma}$ d{\'e}fini par
$ r_{\sigma}:= \sup \{ |\lambda|;  \lambda \in {\rm Sp}({\bf H}_{\sigma}), 
\lambda \not = \lambda(\sigma) \} 
$,
v{\'e}rifie  $r_{\sigma}<\lambda(\sigma)$.

$(3)$ Pour $\sigma\in\Sigma_{0}$, on
 d{\'e}finit la pression par $\Lambda(\sigma = \log\lambda(\sigma)$.
Notons $\Lambda^{''}_{s^{2}},\Lambda^{''}_{\tau^{2}}$ les
d{\'e}riv{\'e}es partielles d'ordre 2 de la fonction $\Lambda(s,i\tau)$.
Au point $(1,0)$ la pression est strictement convexe en $s$, c.{\`a}.d
$\Lambda^{''}_{s^{2}}(1,0)>0$. De plus, si c n'est pas constante,
alors la pression est strictement convexe par rapport {\`a} $i\tau$ en $(1,0)$,
c.{\`a}.d $\Lambda^{''}_{\tau^{2}}(1,0)>0$

$(4)$ {\rm [La fonction  $i\tau \mapsto \sigma(i\tau)$.] } 
Il existe un voisinage complexe $\cal W$ de $0$ et
une unique fonction $\sigma:\mathcal{W} \longrightarrow \C$\; 
tels que\; $\lambda(\sigma(i\tau), i\tau) 
= 1$,
cette fonction est analytique, de plus  
$\sigma(0) = 1, \mbox{et}\;\; \sigma^{''}(0)\neq 0$.  
\end{prop}
\pn Pour la preuve voir (\cite{BVVB2003}, Sec.2).

\medskip
\pn {\bf Condition UNI:}
On dit que le syst{\`e}me dynamique $T$ v{\'e}rifie la condition {\em UNI} si
toute branche inverse de $T$ s'{\'e}tend en une fonction $\cal{C}^{3}$, et
si pour tout $h, k$ deux branches inverses de la m{\^e}me profondeur,  on note 
$$\Psi_{h,k}(x) := \log\frac{|h^{'}(x)|}{|k^{'}(x)|},
\hspace{2cm} \Delta(h,k) :=
\inf_{x\in\cal{I}}\big|\Psi^{'}_{h,k}(x)\big|$$
et pour tout $\eta>0$, 
$$\mathcal{J}(h,k) := \bigcup_{k\in\h^{n},\Delta(h,k)\leq\eta}k(\cal{I}),$$ 
alors,

$(a)$ Pour tout $0<a<1$ on a
$\big|\mathcal{J}(h,\rho^{an})\big|\ll\rho^{an}, \forall n, \forall
h\in\h^{n}$.

\medskip
$(b)\;\;\sup\big\{\big|\Psi^{''}_{h,k}(x)\big| ; n\geq 1,
h,k\in\h^{n}, x\in\cal{I}\big\} <\infty.$

\md
\pn Afin de supprimer l'effet du facteur $|s|$ on d{\'e}finit la norme
suivante :
\[\|f\|_{1,t} :=\|f\|_0+\ds\frac{\|f\|_1}{|t|} = \sup|f|+\ds\frac{\sup|f^{'}|}{|t|}\]

\begin{thm}
$[$Estimations {\`a} la Dolgopyat$]\; ($\cite{BVVB2003} Sec.2$)$.
Soient $({\cal I}, T, c)$ avec $T$ dans la bonne classe et $c$ {\`a}
croissance mod{\'e}r{\'e}e, et $\rho$
le taux de contraction, et telle que la condition {\em UNI} ait lieu.  
Soit ${\bf H}_{s, w}$
son op{\'e}rateur de transfert pond{\'e}r{\'e} agissant sur $\cal C^1(\cal I)$.

\pn Pour tout $r >0$, il existe un voisinage complexe
$\Sigma_1=]1-\alpha, 1+\alpha[$ de $1$ et $M>0$ tels que, 
 pour tout $s=\sigma+it$,  
avec $\sigma\in  \Sigma_1$ et  $|t|\ge 1/ \rho^2$, et tout $\tau\in\R$    
 \[\|(I-{\bf H}_{s, i\tau})^{-1}\|_{1,t} \le M\;|t|^r \, .\]
\end{thm}

Consid{\'e}rons l'ensemble $\Omega _N := \{(u, v)\in\Z^{+}_{\star};\;
gcd(u, v) = 1, u \leq v\leq N\}$, muni de la probabilit{\'e} uniforme $P_N$.
Afin d'{\'e}tudier la distribution du co{\^u}t total $C(u,v)$ associ{\'e} {\`a} un
certain co{\^u}t $c$ ({\`a} croissance mod{\'e}r{\'e}e), on d{\'e}finit sa ``fonction
g{\'e}n{\'e}ratrice des moments" sur  $\Omega _N :$
\[E_N[\exp( i\tau C)]:=\frac{\Phi _{i\tau}(N)}{\Phi _{0}(N)},\]
o{\`u}  $\Phi_{i\tau}(N) = \Phi_{c, i\tau}(N)$  est la valeur cumul{\'e}e de
$\exp(i\tau C)$ sur $\Omega _N$ :
$$\Phi_{i\tau}(N) := \ds \sum_{(u, v)\in\Omega _N}\exp [i\tau C(u,
v)] \;,\;\;\;\;\;\Phi_{0}(N) =  \mid\Omega _N \mid.$$
En suivant le principe d{\'e}fini dans \cite{Vallee1998},
 on peut
remplacer la suite des fonctions g{\'e}n{\'e}ratrices des moments  
 par une {\it s{\'e}rie de Dirichlet,} qu'on appelera 
 {\it la fonction g{\'e}n{\'e}ratrice des  moments de Dirichlet:}

\[S(s, i\tau) := \ds \sum_{(u, v)\in \Omega}
 \frac {1} {v^s}  \exp [i\tau C(u, v)]
 = \ds\sum_{n \ge 1} \frac{c_n(i\tau)}{n^s},\]

\pn o{\`u}  \[\Omega :=\{(u, v)\in\Z^{+}_{\star};\; gcd(u, v) = 1\}\;\mbox{et}\;
c_n(i\tau):=\ds\!\! \sum _{(u, v)\in\Omega_n\, , v=n}\!\!\exp[i\tau C(u, v)].\]
On a
$$\ds \sum_{n\leq N}c_n(i\tau) = \Phi_{i\tau}(N).$$
De plus, il est facile de montrer que ;
\begin{equation}\label{Dirichlet} 
S(2s, i\tau)={\mathbf F}_{s, i\tau} \circ (I- {\mathbf H}_{s, i\tau})^{-1}[1] (0),
\end{equation}

\md
\pn o{\`u} 
\[{\mathbf F}_{s, i\tau} [f](x) := {\mathbf H}_{s, i\tau} [f \cdot 1_{\cup_{ h
      \in \cal F} 
h(\cal I)}](x),\] 

\md
\pn Parmi les r{\'e}sultats obtenus par V.Baladi et B.Vall{\'e}e
(\cite{BVVB2003}, Sec.4) le
{\it CLT} suivant :

\begin{thm}
$[$Th{\'e}or{\`e}me de la limite centrale avec vitesse de conver-gence.$]$ 
Pour les algorithmes  Euclidiens 
$\dji$, $\key$, $\ow$, et tout  co{\^u}t  $c\not\equiv 0$ {\`a} croissance mod{\'e}r{\'e}e,
en posant $\Lambda(\sigma) := \Lambda(\sigma,0)$
la fonction de la proposition 1:
$(a)$ Il  existe $\mu(c)> 0, \delta(c)> 0$ et $K>0$ tels que, pour tout
$N\in\N^{*}$ et tout $y\in\R$
\begin{equation*}
\bigg|P_{N}  \left[(u, v)\;\;| \frac{C(u, v) - \mu (c) \log N}{\delta
  (c)\sqrt{\log N}} \leq y \right]
-\frac{1}{\sqrt{2\pi}} \intop _{- \infty}^{y} e^{-x^{2}/2}
  \;dx\bigg|\leq\!\!\frac{K}{\sqrt{\log N}}. 
\end{equation*}

\pn o{\`u}
$$\mu(c)=2\sigma^{'}(0)\;\;\mbox{et}\;\;\;\delta^{2}=2\sigma^{''}(0).$$
\pn avec $\sigma$ la fonction de la proposition 1.4.
\end{thm}

\md
\pn Pour un intervalle $J$ de $\R$, on note $|J|$ la mesure de
Lebesgue de $J$. Notre r{\'e}sultat montr{\'e} dans la troisi{\`e}me partie, est
le th{\'e}or{\`e}me suivant: 
\begin{thm}$[\mbox{Th{\'e}or{\`e}me de la limite locale.}]$
Pour les algorithmes euclidiens
$\mathcal{G},\;\mathcal{K},\;\mathcal{O}$ et pour toute fonction co{\^u}t
$c$ {\`a} croissance mod{\'e}r{\'e}e, en posant $\mu (c), {\delta
  ^2}(c)$ les constantes du th{\'e}or{\`e}me 2 , on a : $\forall J$ intervalle
de $\R$, $\forall \varepsilon>0,\,\exists N_{0}$ tel que pour tout
$N\geq N_{0}$ et tout $x\in\R$
\begin{equation}\label{local}
\bigg|\!\!\sqrt {\log N} P_{N}\!\bigg[\!(C(u,v)-\mu(c)\log N - \delta(c)x\sqrt{\log
  N}) \in J\bigg]\! - |J| \frac{e^\frac{-x²}{2}}{\delta (c)
  \sqrt {2 \pi}}\!\bigg|\! <\!\varepsilon.  
\end{equation}
\end{thm}

\section {Estimations de la fonction g{\'e}n{\'e}ratrice des moments.}

Parmi les cons{\'e}quences les plus utiles de la proposition 1 est que
pour $(s,i\tau)\in\mathcal{W}_{1}$ un voisinage complexe de $(1,0),$ on a ${\bf H}_{s,i\tau}
= \lambda (s, i\tau){\bf P}_{s,i\tau} + {\bf N}_{s,i\tau}$  o{\`u} ${\bf P}_{s,i\tau}$ est
la projection spectrale associ{\'e}e {\`a} $\lambda (s, i\tau)$ et le rayon
spectral de ${\bf N}_{s,i\tau}$ est $\leq \theta$, avec $r_1<\theta<1$
($r_1$ de la proposition 1).

\pn On montre de plus que pour $(s,i\tau)\in\W_1$ 
\[(I -{\bf H}_{s,  i\tau})^{-1} =   \frac {\lambda(s, i\tau)} {1- \lambda(s, i\tau)}\,  
{\bf P}_{s,   i\tau} +   (I -{\bf N}_{s,  i\tau})^{-1}.\]
 a  pour seule singularit{\'e} dans
$\W_1$  un p{\^o}le  simple en chaque point $(s= \sigma (i\tau), i\tau)$,   
avec r{\'e}sidu, l'op{\'e}rateur non nul
 
\[ {\bf R}(i\tau) :=  \frac {-1} {\lambda'_{s}(\sigma(i\tau), i\tau)} \,    
{\bf P}_{\sigma(i\tau), i\tau}.\]

\md
 On veut exprimer $E_N[\exp( i\tau C)]$ en quasi-puissance, pour cela on introduit
un autre mod{\`e}le probabiliste\; $\big(\overline\Omega_{N}(\xi),
\overline P_{N}(\xi)\big)$ avec $\overline\Omega_{N}(\xi)=\Omega_{N}$: 
on fixe une
fonction\; 
$t \longmapsto \xi (T)$, avec\; $0
\leq\xi (T)\leq 1$, puis pour un entier $N$, on choisit
uniform{\'e}ment un entier $Q$ entre $N-\lfloor N \xi (N) \rfloor$ et
$N$, ensuite on choisit un {\'e}l{\'e}ment $(u, v)$ dans $\Omega _{Q}$.

\md
\pn Dans ce qui suit, la notation $A(l)=O(B(l))$ signifie qu'il existe
$M>0$, tel que pour tout $l$; $|A(l)|\leq M|B(l)|$.   

\md
\pn On note  $\overline E_N[\exp( i\tau C)]$\; la fonction g{\'e}n{\'e}ratrice des moments
d{\'e}finie sur $\big(\overline\Omega _N(\xi),\overline
P_{N}(\xi)\big)$. V.Baladi et B.Vall{\'e}e \cite{BVVB2003} obtiennent : 

\begin{lem}
Consid{\'e}rons l'un des trois Algorithmes $\dji, \key, \ow$. Il existe
$0<\alpha_{0}<1/2=\sigma_0$ (avec $\sigma_0$ de la d{\'e}finition 2)
tel que en posant $\xi(N) = N^{-\alpha_{0}}$ on ait :

$(a)$ La distance entre les distributions $\overline P _{N}(\xi)$ et $P_N$ est
$O(\xi (N))$. 

$(b)$ La fonction g{\'e}n{\'e}ratrice des moments $\overline E_{N}$  de $C$
s'exprime en quasi-puissance. Plus pr{\'e}cis{\'e}ment, il existe 
$<\widehat\alpha_{0}<\alpha_{0}<1/2$
tel que pour toute fonction co{\^u}t {\`a} croissance mod{\'e}r{\'e}e on ait :
$1/2=\sigma_{0}<\widehat\alpha_{0}<\alpha_{0}$ (avec $\sigma_{0}$ de la d{\'e}finition 2)
$$\overline E_{N}[\exp (i\tau  C)] =  \frac { E(i\tau)}{E(0) \sigma (i\tau)}  
N^{2 (\sigma(i\tau)-\sigma(0))} [1 + O (N^{-\widehat\alpha_{0}})],$$
avec un $O-\mbox{terme}$ uniforme par rapport {\`a} $N \longrightarrow \infty$,
  $\tau$ proche de $0$, et  $$E(i\tau) = {\bf F}_{\sigma(i\tau),i\tau}\circ
{\bf R}(i\tau)[1](0).$$
\end{lem}

\pn {\bf Preuve}.(-Esquisse - on r{\'e}f{\`e}re {\`a} \cite{BVVB2003} Sec.4 pour les d{\'e}tails.)
La premi{\`e}re partie du lemme d{\'e}coule de la d{\'e}finition de
$\overline P_{N}(\xi)$ et du fait que $|\Omega_{N}|= KN^{2}(1+O(\log
N/N))$, avec $K>0$ bien d{\'e}fini pour chacun des trois algorithmes,
(voir \cite{BVVB2003} Sec.4.4).

\md
\pn Posons 
\[\Psi_{i\tau}(T) := \ds\sum_{n\leq T} c_{n}(i\tau)(T-n) = \ds\sum_{N\leq
  T}\ds\sum_{n\leq N}c_{n}(i\tau) = \ds\sum_{N\leq
  T}\Phi_{i\tau}(N).\]

\pn En appliquant le th{\'e}or{\`e}me de Cauchy sur la s{\'e}rie de Dirichlet $S(s, i\tau)$
et le rectangle \[U(i\tau) = \{s \ ; \ \Re s = 1
\pm\widehat\alpha\}\times\{s \ ; \ 
\Im s = \pm U\},\] avec $\widehat\alpha_{0}<\widehat\alpha$,
($s\longmapsto S(s,i\tau)$\; 
{\'e}tant m{\'e}romorphe sur
  $U(i\tau)$ lorsque $\tau$ est proche de $0$),
 puis la formule de Perron d'ordre 2 qui transforme
l'int{\'e}grale sur le rectangle $U(i\tau)$ en une int{\'e}grale sur une droite
verticale, on obtient la
formule de quasi-puissance pour la s{\'e}rie $\Psi_{i\tau}(T)$. Ensuite la
relation ;

\begin{equation*}
\begin{split}
\overline\Phi_{i\tau}(N) :&=
\ds\frac{1}{N\lfloor\xi(N)\rfloor}\ds\sum_{Q=N-\lfloor N\xi(N)\rfloor}^{N}\ds\sum_{n\leq
Q}c_{n}(i\tau)\\
&=\ds\frac{1}{N\lfloor\xi(N)\rfloor}\bigg[\Psi_{i\tau}(N)-\Psi_{i\tau}\big(N-\lfloor\xi(N)\rfloor
\big)\bigg],
\end{split}
\end{equation*}
\pn nous permet de transmettre la quasi-puissance {\`a} $\overline\Phi_{i\tau}(N)$,
puis {\`a} $\oen$.\;$\Box$

\begin{defn}
Une fonction $c$ est dite {\it r{\'e}seau} si elle est non nulle et s'il
existe $L_c, L_0>0$ tels
que $L_0/L_c$ soit irrationnel, et $(c-L_0)/L_c$ {\`a} valeurs enti{\`e}res. Le plus grand de ces $L_c$ est appel{\'e}
largeur de $c$.
\end{defn}
\pn Le lemme suivant est une petite g{\'e}n{\'e}ralisation du lemme 15 de \cite{BVVB2003}.
\begin{lem}
On consid{\`e}re l'un des algorithmes $\dji, \key, \ow$. Pour toute
 fonction co{\^u}t
 $c$ {\`a} croissance mod{\'e}r{\'e}e, pour
tout $0<L<\infty$ $($dans le cas o{\`u} $c$ 
est une
fonction r{\'e}seau de largeur $L_{c}$, on prend $0<L\leq\pi/L_{c})$,
et tout $0<\widetilde{\nu}_0<L$,\; il existe $\gamma_{0}=\gamma_0(L, \widetilde{\nu}_0)>0$,
$Q= Q(L, \widetilde{\nu}_0)>0$ 
tels que pour tout\;$|\tau|\in[\widetilde{\nu}_0, L]$ on ait pour 
$\xi(N)=N^{-\widehat\alpha_0}$  
$$\overline{E}_{N}[exp (i \tau  (C(u, v))]\leq\;Q\;N^{-\gamma_0},
\hspace{1cm}\forall N\in \N.$$
\end{lem}

\md
\pn {\bf Preuve}. 
Soit $r>0$, le th{\'e}or{\`e}me 1 assure l'existence de $\alpha>0$ tel que
pour tout $s$ avec $\Re s=\sigma\leq|1-\alpha|$ et $|\Im s|\geq 1/\rho^{2}$ et pour
$\tau$ arbitraire ;

\begin{equation*}
\|\big(I- {\bf H}_{s,i\tau}\big)^{-1}\|_{1,t}\leq
M\;|t|^{r}.
\end{equation*} 

\md
\pn Supposons que $|t|\leq1/\rho^{2}$ et $\tau\in[\widetilde{\nu}_{0}, L]$.
La proposition 1.(1), 1.(3) et la condition UNI impliquent que
$1\notin Sp{\bf H}_{1+it,i\tau}$ (voir prop.1 de \cite{BVVB2003}). Donc
d'apr{\`e}s la th{\'e}orie de la perturbation
de parties finies du spectre il existe\;
$0<\gamma_{1}<\alpha$ \;tel que sur l'ensemble compact 
$$\{(s,\tau)\in\C\times\R;\hspace{0.3cm}|\sigma-1|\leq\gamma_{1},|t|\leq
1/\rho^{2},\;|\tau|\in[\tilde\nu_{0}, L]\},$$
\vskip2pt
\pn on ait $1\notin Sp{\bf H}_{\sigma+it,i\tau}$. En effet la 
fonction $s\longmapsto {\bf H}_{s,i\tau}$ est
analytique sur cet ensemble. D'o{\`u} l'existence de
$\widetilde{Q} = \widetilde Q(\tilde\nu_{0}, L)$ tel que :
\begin{equation*}
\|\big(I- {\bf H}_{1\pm\gamma_{1}+it,i\tau}\big)^{-1}\|_{1,t}\leq\widetilde Q.
\end{equation*} 

\md
\pn Par cons{\'e}quent, pour tout $|\tau|\in[\widetilde\nu_{0}, L]$, il existe
$Q=Q(\widetilde{\nu}_{0},L)$ tel que ;
\begin{equation}\label{Dolg}
\|\big(I- {\bf H}_{1\pm\gamma_{1}+it,i\tau}\big)^{-1}\|_{1,t}\leq
Q\max(1,|t|^{r}),\;\;\;\forall t\in\R.
\end{equation} 

\md
\pn Gr{\^a}ce {\`a} (\ref{Dirichlet}), on transforme (\ref{Dolg}) en une estimation de $S(s,i\tau)$ qui, en
tant que fonction de $s$, est analytique sur le rectangle $\widetilde U(i\tau)
= \{s;\;\Re s = 1\pm\gamma_{1}\}\times\{s;\;\;\Im s = \pm U\}$ (avec
$0<\gamma_{1}<\widehat\alpha, U >0$).

\md
\pn Le th{\'e}or{\`e}me de Cauchy et la formule de Perron nous
permettent de d{\'e}duire la d{\'e}croissance de $\overline\Phi_{i\tau}(N)$ et par
cons{\'e}quent celle de $\oen$ (voir\cite{BVVB2003} Sec.5).\;$\Box$  
 
\section {Preuve du th{\'e}or{\`e}me 3.}

\pn Posons $n=\log N$\; et\; $q_{x}(n) = \mu(c)n - \delta(c)x\sqrt{n}$.

\pn Rappelons que, par le lemme 1, $\big|\overline P_N(\xi) - P_N\big| =
O(e^{-n\alpha_{0}})$, il suffit alors de d{\'e}montrer (\ref{local}) pour $\overline P_{N}$.  

\pn Soient $\overline m_{n}=\overline m_{x,n}$ une suite de mesures d{\'e}finies
sur la tribu des Bor{\'e}liens de $\R$ par :
\begin{equation*}
\overline m_{n}(J) :=\overline P_{N}\bigg[\big(C(u,v)-q_{x}(n)\big)\in J\bigg], 
\end{equation*}

\pn et $m=m_x$ la mesure d{\'e}finie par :
\begin{equation*}
m(J) :=\ds\frac{e^{-x^{2}/2}}{\delta(c)\sqrt{2\pi}}|J|.
\end{equation*}

\pn On suit la m{\'e}thode de Breiman (\cite{Brei}, Chp.10.2) : Pour montrer que 
$\ds\sqrt n\;\overline m_n{\mathop{\longrightarrow}^{w}\limits} m$, il suffit
de montrer que pour toute fonction $\psi$ non-n{\'e}gative continue et
dont la transform{\'e}e de Fourier $\widehat\psi$ est {\`a} support compact on ait :
\begin{equation*}
\sqrt n\intop\psi\;d\overline m_{n}\longrightarrow\intop\psi\;dm.
\end{equation*}

\pn Ce qui signifie que pour tout 
$\varepsilon$ fix{\'e}, il existe $n_0\in\N$ ($n_0$ ind{\'e}pendant de $x$) tel
que pour tout $N\geq N_0$ ;
\begin{equation*}
\bigg|\sqrt{n}\;\oen\big[\psi\big(C(u,v)-q_{x}(n)\big)\big] -
\frac{e^{-x^{2}/2}}{\sqrt{2\pi}\delta(c)}\intop\psi(y)dy\bigg|<\varepsilon,
\quad\forall x\in\R.
\end{equation*}

\pn Soit $[-L, +L]$ contenant le support de $\widehat\psi$. On a  
\begin{equation*}
\begin{split}
\sqrt{n}\;\oen\big[\psi\big(C(u,v)&-q_{x}(n)\big)\big]\\
&=\frac{\sqrt{n}}
{2\pi}\intop_{-L}^{+L}\widehat\psi(\tau)
\oen\big[\exp\big(i\tau\big(C(u,v)-q_{x}(n)\big)\big)\big]\;d\tau\\
&=\frac{\sqrt{n}}{2\pi}\!\intop_{-L}^{+L}\widehat\psi(\tau)\;e^{-i\tau
q_{x}(n)}\oen\big[\exp\big(i\tau\big(C(u,v)\big)\big]\;d\tau\\
&=:I^{(n)}.
\end{split}
\end{equation*}

\pn Soit $0<\nu_0<\widetilde{\nu}_0$ avec $\widetilde\nu_0$ assez
petit (comme dans le lemme 2). 
D{\'e}composons l'intervalle $[-L, +L]$ en
$|\tau|\leq\nu_0$\;et\;$|\tau|\in[\nu_0, L]$, ainsi $I^{(n)}$ se
d{\'e}compose en $I_0^{(n)} + I_1^{(n)}.$

\pn Montrons d'abord que $I^{(n)}_{1} \longrightarrow 0$ 
 lorsque $n\longrightarrow\infty$.

\md
\pn En appliquant le lemme 2, on obtient (rappelons que
$\gamma_{0}=\gamma_{0}(L, \widetilde\nu_{0})$)

\begin{equation*}
\begin{split}
|I^{(n)}_{1}| &\leq\frac{\sqrt n}{2\pi}\;Q\;N^{-\gamma_0} \intop
_{|\tau|\in[\nu_0, L]} 
\big|\widehat\psi(\tau)\big|\;d\tau\\
&\leq\frac{Q}{2\pi} \sqrt n\;e^{-n\gamma_{0}}\intop
_{|\tau|\in[\nu_0, L]} 
\big|\widehat\psi(\tau)\big|\;d\tau\\
&\leq\frac{\widetilde{K}}{n} \intop
_{|\tau|\in[\nu_0, L]} 
\big|\widehat\psi(\tau)\big|\;d\tau =\frac{\tilde{C}}{n}<\varepsilon/2.\\
\end{split}
\end{equation*} 

\pn Il suffit de prendre
 $N > \exp\left[\widetilde{C}\varepsilon^{-1}\right]$, \;o{\`u}\; $\widetilde
 C=\tilde C (L,\;\widetilde\nu,\;\psi)$, $\tilde C$ ind{\'e}pendant de $x$.

\pn Calculons  $I^{(n)}_{0}$. 
On suit la m{\'e}thode \cite{BVVB2003} (sec.5).
      
\pn Rappelons d'abord que d'apr{\`e}s la proposition 1.4,
$\sigma^{''}(0)\neq 0$,
d'o{\`u} pour $\nu_{0}$ suffisamment petit, on a $\delta_0 :=
 \inf\{|\Re\sigma^{''}(\tau)|;\;\tau\in[-\nu_{0}, \nu_{0}]\}>0$.

\pn Posons $ \tau_n := \left(\ds\frac{\log n}{\delta_{0} n}\right)^{1/2}$ et
 d{\'e}composons l'intervalle $[- \nu_0, +\nu_0]$  en $|\tau| \leq \tau _n$
 et $|\tau| \in [\tau_n , \nu_0]$.

\md
\begin{equation*}
\begin{split}
I^{(n)}_{0} &=\frac{\sqrt n}{2\pi}\;\intop _{|\tau|\leq\tau _n} 
\widehat { \psi}(\tau)\;
e^{-i\tau  q _ {x} (n)}\overline{E}_{N}[exp(i \tau C(u, v))]\;d \tau\\
&\;\;\;\;\;\;\;\;\;\;+\frac{\sqrt n}{2\pi} \intop _ {|\tau|\in [\tau _n, \nu
  _0]} \widehat {\psi}(\tau)\;e^{-i\tau  q _ {x} (n)}
\overline{E}_{N}[exp(i \tau C(u, v))]\;d \tau.\\
\end{split}
\end{equation*} 
Le 2{\`e}me terme est {\'e}gal {\`a} $O\big(1/\sqrt n\big)$, en effet, en
rappelant l'expression de $\overline E_{N}$ dans le lemme 1, et que la
fonction $g : z\longmapsto \sigma(z)  - 1 - z\sigma^{''}(0)$ admet un
point col en $z=0\;\big(g^{'}(0) = 0\;\; 
\mbox{et}\; g^{''}(0)\neq 0\big),$ on a :   

\[\ds
 \begin{array}{l}
\ds\frac{\sqrt n}{2\pi}\;\bigg|\!\!\!\! \intop_{|\tau|\in [\tau _n, \nu _0]} 
\wpsi(\tau)\;e^{-i\tau  q_{x}(n)}
\overline{E}_{N}[exp(i\tau C(u, v))]\;d \tau\bigg| \\\\
\ds=\frac{\sqrt n}{2\pi} \bigg|\intop_{|\tau|\in [\tau _n, \nu _0]}\widehat
{\psi}(\tau)\;e^{-i\tau \delta (c)x\sqrt n} 
N^{-i\tau\mu (c)}\overline{E}_{N}[exp(i \tau C(u, v))]\;d \tau \bigg|\\\\
\ds\leq \frac{\sqrt n}{2\pi} \!\sup_{|\tau|\in [\tau _n, \nu _0]}\! \bigg| N^{-i \tau
\mu (c)}\overline{E}_{N}[exp(i\tau
C(u, v))]\!\bigg|\;\; \bigg|\!\!\!\!\intop_{|\tau|\in [\tau _n,
  \nu _0]}\!\!\!
\!\!\!\widehat
{\psi}(\tau)\;e^{-i\tau \delta (c)x\sqrt n}\;d \tau\bigg|\\\\
\end{array}
\]
\[
\ds\leq\frac{\sqrt n}{2\pi}\!\!\! \sup _ {|\tau|\in [\tau _n, \nu _0]}\! \bigg|
e^{2n(\sigma(i\tau)-1-i\tau\sigma^{'}(0))}\frac{E(i\tau)}{E(0)\;
\sigma(i\tau)}(1+O(e^{-\widehat\alpha_{0}
n}))\bigg|\! 
\intop_{\tau\in\R}\!\!|\widehat\psi(\tau)|\;d\tau.\]   

\pn On v{\'e}rifie ais{\'e}ment que les fonctions : $$f^{(n)}_{1}
:\tau\longmapsto\exp\big[2n\big(\sigma(i\tau)-1-i\tau\sigma^{'}(0)\big)\big]$$
et
$$f^{(n)}_{2}
: \tau\longmapsto\frac{E(i\tau)}{E(0)\sigma(i\tau)}
\big(1+O(e^{-\widehat\alpha_{0}n})\big)$$
v{\'e}rifient pour $|\tau|\leq\nu_0$ : 
\begin{equation}\label{Decrois}
\big|f^{(n)}_{1}(\tau)\big| = O\big(e^{-n\tau^{2}\delta_{0}}\big)
\end{equation}
et
\begin{equation}\label{Taylor}
\big|f^{(n)}_{2}(\tau)\big|= O\big(1+|\tau|+
e^{-n\widehat\alpha_{0}}\big).
\end{equation}

\pn D'o{\`u}
\begin{equation*}
\begin{split}
\frac{\sqrt n}{2\pi}\;\bigg|\!\!\!\! &\intop_{|\tau|\in [\tau _n, \nu _0]} 
\wpsi(\tau)\;e^{-i\tau  q_{x}(n)}
\overline{E}_{N}[exp(i\tau C(u, v))]\;d \tau\bigg| \\\\
&\leq\widetilde{D}(\psi)\frac{\sqrt n}{2\pi}\sup_{|\tau|\in [\tau _n, \nu _0]} 
\bigg| e^{-n\tau^{2}\delta_{0}}(1+|\tau|+e^{-n\widehat\alpha_{0}})\bigg|\\\\
&\leq\sqrt n\; M\; e^{-n\tau^{2}_{n}\delta_{0}}\\\\
&=\sqrt n\;\frac{M}{n} = O\big(\ds1/\sqrt n\big).\\
\end{split}
\end{equation*}

\pn Avec $M$ d{\'e}pendant de $\psi,\;\widetilde\nu_0$ et ind{\'e}pendant de $x.$

\md
\pn Rappelons que, par le th{\'e}or{\`e}me 2, ${\mu(c)\!=\!2\sigma^{'}(0)\;\mbox{et}
\;\delta^{2}(c)\!=\!2\sigma^{''}(0)}$; ainsi,

\begin{equation*}
\begin{split}
&\ds \widetilde{I}^{(n)}_{0}:=\frac{\sqrt n}{2\pi}
\intop_{|\tau|\leq\tau_{n}}\!\wpsi(\tau)\;e^{-i\tau
  q_{x}(n)}\overline{E}_{N}[exp(i\tau C(u, v))]\;d \tau \\\\
&=\!\frac{\sqrt
  n}{2\pi}\!\!\!\!\!\!\intop_{|\tau|\leq\tau_{n}}\!\!\!\!\!\wpsi(\tau)
e^{-i\tau
  q_{x}(n)}e^{2n(\sigma(i\tau)-1)}
\frac{E(i\tau)}{E(0)\sigma(i\tau)} 
\big(1+O(e^{-\widehat\alpha_{0}n})\big)d\tau\\\\
&=\frac{\sqrt n}{2\pi}\!\!\!\!\intop_{|\tau|\leq\tau_{n}}
\!\!\!\!\!\!\wpsi(\tau)\;e^{-i\tau\delta(c)x\sqrt
n}\;e^{2n(\sigma(i\tau)-1-i\tau\sigma^{'}(0))}
 \ds\frac{E(i\tau)}{E(0)\sigma(i\tau)}
\big(1\!+\!\!O(e^{-\widehat\alpha_{0}n})\big)\;d\tau.
\end{split}
\end{equation*}
\pn Posons 
$$J^{(n)}_{0}:=\frac{\sqrt n}{2\pi}\intop_{|\tau|\leq\tau_{n}}
\wpsi(\tau)\;e^{-i\tau\delta(c)x\sqrt
n}\;e^{2n(\sigma(i\tau)-1-i\tau\sigma^{'}(0))} 
d \tau.$$
\pn Gr{\^a}ce {\`a} (\ref{Decrois}) et (\ref{Taylor}), on a
\begin{equation*}
\begin{split}
\frac{1}{2\pi}\big|\widetilde{I}_{0}^{(n)} - J_{0}^{(n)}\big|&\leq
\frac{1}{2\pi}
\intop_{|\tau|\leq\tau_n}\!\!\!\!\big|\wpsi(\tau)\big||\tau|e^{-n\tau^{2}\delta_0}
+ e^{-n\widehat\alpha_0}\intop_{|\tau|\leq\tau_n}\!\!\big|\wpsi(\tau)\big|
e^{-n\tau^{2}\delta_0} d\tau\\
&\leq\frac{1}{2\pi}\intop_{|\tau|\leq\tau_n}\big|\wpsi(\tau)\big||\tau|\;d\tau
+
\frac{1}{n\widehat\alpha_0}\intop_{|\tau|\leq\tau_n}\big|\wpsi(\tau)\big|\;d\tau\\
&\leq \frac{K_1}{\sqrt n} + \frac{K_2}{n} = O\big(1/\sqrt n\big). 
\end{split}
\end{equation*}
O{\`u} $K_{1}, K_{2}$ d{\'e}pendent de $\psi.$

\pn Ainsi,
\begin{equation*}
\widetilde{I}_{0}^{(n)} =\frac{\sqrt
  n}{2\pi}\!\!\intop_{|\upsilon|\leq\tau_n\sqrt
  n}\!\!\wpsi(\frac{\upsilon}{\sqrt n})
\;e^{-i\tau\delta(c)x\sqrt
  n}\;e^{2n(\sigma(i\frac{\upsilon}{\sqrt
  n})-1-i\frac{\upsilon}{\sqrt n}\sigma^{'}(0))}d\upsilon +O\big(1/\sqrt n\big). 
\end{equation*}

\pn D'autre part, 
\[\bigg|e^{2n(\sigma(i\frac{\upsilon}{\sqrt
  n})-1-i\frac{\upsilon}{\sqrt n}\sigma^{'}(0))} -
e^{\frac{\delta^{2}(c)\upsilon^{2}}{2}}\bigg|=
  O\bigg(\ds\frac{\upsilon^3}{\sqrt n}\bigg)\]

\pn d'o{\`u} par le th{\'e}or{\`e}me de la convergence domin{\'e}e de Lebesgue, on a
pour tout $n$  suffisamment grand et ind{\'e}pendant de $x$ :

$$
\bigg|\widetilde{I}^{(n)}_{0} - \wpsi(0)\intop_{\R}
e^{-i\delta(c)x\upsilon}e^{-\delta^{2}(c)
\upsilon^{2}/2}d\upsilon\bigg|<\varepsilon/2,\hspace{1cm}\forall x\in\R.$$

\pn Ce qui signifie que pour tout $\varepsilon$ fix{\'e}
et $n$ suffisamment grand, 

$$\bigg|I^{(n)}_{0} -  
\frac{e^{-x^{2}/2}}{\delta(c)\sqrt{2\pi}}\intop\psi(y)dy\bigg|<\varepsilon/2,$$
uniform{\'e}ment en $x.\;\Box$

\md
\pn {\bf Remerciements.} Je remercie vivement V.Baladi et B.Vall{\'e}e pour les
lectures de cet article et leurs remarques enrichissantes, S.Gou{\"e}zel
pour sa remarque sur la m{\'e}thode de Breiman, S.Khemira et L.Pharamond pour
leur pr{\'e}cieuse aide informatique. 

\md
\goodbreak

\end{document}